\documentclass[bibother]{asl}
\usepackage{rotate}

\usepackage[UKenglish]{babel}
\usepackage[utf8]{inputenc}
\usepackage{verbatim}

\usepackage{theoremref}
\usepackage[]{hyperref}
\usepackage{graphicx}

\usepackage{tikz-cd} 

\usepackage{dsfont}
\usepackage{enumitem}

\usepackage{thmtools}

\usepackage{amssymb}
\usepackage{relsize}

\makeatletter

\def\Autoref#1{%
	\begingroup
	\edef\reserved@a{\cpttrimspaces{#1}}%
	\ifcsndefTF{r@#1}{%
		\xaftercsname{\expandafter\testreftype\@fourthoffive}
		{r@\reserved@a}.\\{#1}%
	}{%
		\ref{#1}%
	}%
	\endgroup
}
\def\testreftype#1.#2\\#3{%
	\ifcsndefTF{#1autorefname}{%
		\def\reserved@a##1##2\@nil{%
			\uppercase{\def\ref@name{##1}}%
			\csn@edef{#1autorefname}{\ref@name##2}%
			\autoref{#3}%
		}%
		\reserved@a#1\@nil
	}{%
		\autoref{#3}%
	}%
}
\makeatother

\theoremstyle{plain}
\newtheorem{theorem}{Theorem}[section]
\newtheorem{corollary}[theorem]{Corollary}
\newtheorem{lemma}[theorem]{Lem\-ma}
\newtheorem{proposition}[theorem]{Prop\-o\-si\-tion}

\newtheorem{question}[theorem]{Question}

\newtheorem{fact}[theorem]{Fact}

\theoremstyle{definition}
\newtheorem{definition}[theorem]{Definition}

\newtheorem{remark}[theorem]{Remark}
\newtheorem{example}[theorem]{Example}
\newtheorem{notation}[theorem]{Notation}

\newcommand{\N}[0]{\mathbb{N}}
\newcommand{\Z}[0]{\mathbb{Z}}
\newcommand{\Q}[0]{\mathbb{Q}}
\newcommand{\R}[0]{\mathbb{R}}

\renewcommand{\AA}[0]{\mathcal{A}}
\newcommand{\II}[0]{\mathcal{I}}
\newcommand{\OO}[0]{\mathcal{O}}
\newcommand{\MM}[0]{\mathcal{M}}
\newcommand{\NN}[0]{\mathcal{N}}
\newcommand{\supp}[0]{\mathrm{supp}}

\renewcommand{\L}[0]{\mathcal{L}}
\newcommand{\Lor}[0]{\mathcal{L}_{\mathrm{or}}}
\newcommand{\Log}[0]{\mathcal{L}_{\mathrm{og}}}

\newcommand{\Lsp}[0]{\mathcal{L}_{\mathrm{sp}}}

\newcommand{\ul}[1]{\underline{#1}}
\newcommand\restr[2]{{
		\left.\kern-\nulldelimiterspace 
		#1
		\vphantom{\big|} 
		\right|_{#2}
}}
\newcommand{\SPn}[0]{\mathrm{SP}_n}
\newcommand{\vmin}[0]{{v_{\min}}}
\newcommand{\vnat}[0]{{v_{\mathrm{nat}}}}

\DeclareMathOperator{\HProd}{{\mathlarger{\mathlarger{\mathbf{H}}}}}

\DeclareMathOperator{\rk}{rk}

\DeclareMathOperator {\drk}{drk}

\newcommand\blfootnote[1]{
	\begingroup
	\begin{NoHyper}
	\renewcommand\thefootnote{}\footnote{#1}
	\addtocounter{footnote}{-1}
	\end{NoHyper}
	\endgroup
}

\usepackage{xcolor}

\def\confname{}
\def\copyrightyear{2000}

\title{Definable ranks}

\author[L.~S.~Krapp]{Lothar Sebastian Krapp}
\revauthor{Krapp, Lothar Sebastian}
\author[S.~Kuhlmann]{Salma Kuhlmann}
\revauthor{Kuhlmann, Salma}
\author[L.~Vogel]{Lasse Vogel}
\revauthor{Vogel, Lasse}

\address{Institut f\"{u}r Interdisziplin\"{a}re Sprachevolutions\-wissenschaft, Uni\-ver\-si\-tät Zürich, 8050 Zürich, Switzerland \& Fachbereich Mathematik und Sta\-tistik, Universität Konstanz, 78457 Konstanz, Germany}
\email{sebastian.krapp2@uzh.ch}
\urladdr{https://www.isle.uzh.ch/en/staff/krapp.html}

\address{Fachbereich Mathematik und Statistik, Universität Konstanz, 78457 Konstanz, Germany}
\email{salma.kuhlmann@uni-konstanz.de}
\urladdr{https://www.mathematik.uni-konstanz.de/kuhlmann/}

\address{Fachbereich Mathematik und Statistik, Universität Konstanz, 78457 Konstanz, Germany}
\email{lasse.vogel@uni-konstanz.de}
\urladdr{https://www.math.uni-konstanz.de/\urltilde vogel/}

\begin{document}
	
\begin{abstract}
	We introduce the notion of the definable rank of an ordered field, ordered abelian group and ordered set, respectively. We study the relation between the definable rank of an ordered field and the definable rank of the value group of its natural valuation. Similarly, we compare the definable rank of an ordered abelian group to that of its value set with respect to the natural valuation. We describe the definable rank on the group-level by characterizing the definable convex subgroups. We also give a detailed comparison of field- and group-level, in particular for ordered fields with henselian natural valuation. We investigate definability of final segments in ordered sets and introduce definable condensation as a tool for further study.
\end{abstract}

\maketitle

\section{Introduction}
	Consider the rank of a (linearly) ordered field as the collection of its non-trivial convex valuation rings (cf.\ \cite[page~93]{PriessCrampe83}) and the rank of an ordered abelian group as the collection of all proper convex subgroups (cf.\ \cite[page 56]{Fuchs63},\cite{PriessCrampe83},\cite[page~26]{engler}). These collections are ordered by set inclusion, enabling us to consider ranks as ordered sets.
	Classically, their order types are used as tools to capture essential valuation-theoretic information from the structures under consideration (cf.\ \cite{ZariskiSamuel58}), we will however need to use the exact sets for our purposes.
	More recently, the notion of rank has been considered for more general  ordered algebraic structures, where the specific notion of rank depends on the algebraic complexity of the structure.
	For instance, exponential ranks are crucial in the study of ordered exponential fields (see \cite{kuhlmann}) or ordered transexponential fields (see \cite{KrappKuhlmann25}). Likewise, differential ranks are utilised for the examination of differential fields (see \cite{KuhlmannLehericy19}) or difference ranks for difference fields (see \cite{KuhlmannMatusinskiPoint17}).\blfootnote{Math Subject Classification (2020): Primary 03C64 03C40; Secondary 12J15 12L12 12J10 06F20 06A05. Keywords: real field, definable valuation, definable convex subgroup, definable final segment, henselian.}

	Crucially, there is an order-preserving one-to-one correspondence between the rank of an ordered field and the rank of its value group under the natural valuation (cf.\ \cite[Lemma 3.4]{kuhlmann}). Likewise, there is an order-preserving one-to-one correspondence between the rank of an ordered abelian group and the proper final segments of its value set under the natural valuation, ordered by inclusion (cf.\ \cite[Lemma 3.5]{kuhlmann}). As a result, we define the rank of an ordered set as the collection of its proper final segments.
	
	From a model theoretic perspective, substructures are of particular interest if they are first-order definable. Our aim in this work is to adapt the notion of ranks to the search for distinguished definable substructures. Indeed, we shall establish the \emph{definable} rank as the subset of the rank which only consists of definable substructures.
	
	One of our aims is to study the definable ranks on field-, group- and set-level individually.
	Studying the definable rank of an ordered field is in essence an examination of its definable convex valuations. Of particular note is the relation to ``Shelah's conjecture'' (see \cite[Conjecture 5.34~(c)]{shelah}), which suggests that every strongly dependent/NIP field should carry a definable henselian valuation. In recent time, the pursuit of this conjecture has been a motivation factor for many works on the subject of definable valuations, such as \cite{DittmannJahnkeKrappKuhlmann23, dupont, fehm, krapp, krapp2, KrappKuhlmannLink23, KrappKuhlmannVogel24}. Further work on definable valuations and related considerations, in particular regarding definable henselian valuations, can be found in \cite{delon, FehmJahnke15, hong, jahnke4, jahnke5, prestel}.
	On the group-level, we recompile results from \cite{Schmitt1982} and \cite{delon}. With that we close any remaining gaps to obtain a full characterisation of definable convex subgroups in \autoref{theorem:CharacterizationOfDefinableConvexSubgroups} in the language of \cite{Schmitt1982}. Afterwards we translate most cases into algebraic conditions.
	We also compile initial thoughts and approaches on the matter of definable final segments of ordered sets. We develop the definable condensation in \autoref{def:dfblCondense}, which can be employed as a tool for further study. 

	Our second aim is to study the definable ranks on the three levels in relation to each other.
	An evident question would be, whether we also obtain one-to-one correspondences when we now consider the definable ranks. This is not the case, we show that in general there does not exist a one-to-one correspondence of the definable rank of an ordered field, the definable rank of its value group under the natural valuation and the definable rank of this groups value set under its natural valuation respectively by giving a counterexample with \autoref{xmpl:allRanksDifferent}. We further investigate the relations under the isomorphisms of the full ranks and illustrate in \autoref{xmpl:drkNonIsomorphic}, that even in cases, where the definable ranks are isomorphic as ordered sets, they may not be in correspondence via the restriction of these isomorphisms.

	Despite this, it turns out that group- and field-level are still closely related. Investigation of their connection leads to the first main result \autoref{theorem:mostDfblSubGpCorrespDfblVal}. In particular, if we consider ordered fields with henselian natural valuation, the added tameness of henselian fields lets us refine this result to a complete description of the relation between the two levels, depending on algebraic properties of the field and group. This is elaborated as the second main result in \autoref{theorem:MainTheorem}.\\

	The structure of this works is as follows. In \autoref{section:prelims} we gather necessary preliminaries.
	In \autoref{section:defRanks}, we introduce definable ranks and illustrate their initial disparities.

	In \autoref{section:defCSubgroups}, we investigate the definable rank of an ordered abelian group  and obtain several necessary or sufficient conditions for the definability of a given convex subgroup. 
	This section makes heavy use of the tools from \cite{Schmitt1982}, along with additional insights from \cite{delon}. We use these tools to finish a complete characterisation of definable convex subgroups in terms of Schmitt's machinery in \autoref{theorem:CharacterizationOfDefinableConvexSubgroups} by proving the converse direction of \cite[Corollary 4.2]{delon}. We furthermore extract algebraic conditions for many cases that allow a quicker verification than using the heavy machinery in several cases one might want to consider.

	In \autoref{section:FieldAndGroup}, we use known sufficient conditions for the definability of convex valuation rings from \cite{DittmannJahnkeKrappKuhlmann23} to ultimately compare the definable ranks on field- and group-level. We see through \autoref{theorem:mostDfblSubGpCorrespDfblVal} that, despite the earlier examples, they are in general more closely connected than one might assume. In particular, if we assume that the natural valuation of an ordered field is henselian, then we establish in \autoref{theorem:MainTheorem} that both definable ranks are almost the same and completely characterize when and how they may differ.

	In \autoref{section:finalSegments}, we first study the definability of final segments of an ordered set. In the situation where the ordered set is either dense or discrete, we	give a full characterisation of its definable final segments in \autoref{lem:denseDiscreteSimple}. We then illustrate that all open cases stem from orderings that are locally neither dense nor discrete with \autoref{fact:Rubin}. An approach to this problem would be to employ condensations, i.e.\ quotients of linear orderings that aim to gradually eliminate non-dense parts. To stay compatible under the lens of definability, we develop the definable condensation in \autoref{def:dfblCondense}.\\

	\noindent\textbf{Acknowledgement:} This work is part of the third author's doctoral research project `Dependent ordered structures',	funded by Evangelisches Studienwerk Villigst. It was initiated during the \emph{Tame geometry} conference in February 2025, during which all three authors were generously hosted by CIRM, Marseille. The first author received partial project funding by Vector Stiftung. The authors also want to thank Franz Viktor Kuhlmann for pointing out an error in an earlier version.

\section{Preliminaries}\label{section:prelims}
	We mainly follow the notation and terminology of \cite{Marker02} and \cite{kuhlmann} when applicable. 
	\subsection{Ordered structures}
	We write $\N$ for the set of natural numbers without $0$. We denote the language of ordered rings $\{ +, -, \cdot, 0, 1, < \}$ by $\L_{\mathrm{or}}$, the language of ordered groups $\{ +, -, 0, < \}$ by $\L_{\mathrm{og}}$ and the language of ordered sets $\{ < \}$ by $\L_{\mathrm{<}}$. If no confusion is likely to arise, we omit the symbols $+, -, \cdot, 0, 1$ and simply write $(G,<)$ for the $\L_{\mathrm{og}}$-structure of an ordered abelian group and $(K,<)$ for the $\L_{\mathrm{or}}$-structure of an ordered field. If furthermore the considered ordering is clear, then we further omit the symbol $<$ and only write the domain for the respective structure. When we work with families of structures, we sometimes abuse terminology and omit the subscript on the symbols of the given structures.
	Let $\L$ be a (first-order) language. For an $\L$-structure $\MM$ with domain $M$, we call a subset of $M$ \textbf{definable} if it is $\L$-definable with parameters from $M$; we call it \textbf{$\emptyset$-definable} if it is definable without parameters.
	
	For the entire work, whenever we speak of an ordering we mean a linear ordering. For an ordered set $(I,<_I)$, we say that a subset $J \subseteq I$ is \textbf{convex} (in $I$) if for any $a,b \in J$ and any $c \in I$ with $a <_I c <_I b$ we already have $c \in J$. We say that $J \subseteq I$ is a \textbf{final segment} (of $I$) if for any $a \in J$ and any $ b \in I$ with $a <_I b$, also $b \in J$. We call $I$ \textbf{well-ordered} if every subset of $I$ has a minimum. We call $I$ \textbf{reversed well-ordered} if every subset has a maximum.
	For an ordered set $(I,<)$ and $a, b \in I$ we say $b$ is the \textbf{successor} of $a$ and $a$ is the \textbf{predecessor} of $b$ if $a<b$ and there is no $c \in I$ with $a<c<b$.	
	An ordered set $(I,<)$ is \textbf{dense} if no $i \in I$ has a successor. We call $(I,<)$ \textbf{discrete} if for any $a, b \in I$ with $a<b$ there exists a successor of $a$ and a predecessor of $b$.	
	For any $\gamma \in I$, we define the final segments $\gamma_-:=\{i\in I\mid \gamma\leq i\}$ and  $\gamma_+:=\{i\in I\mid \gamma < i\}$.

	We denote by $(\omega,<)$ the ordered set of non-negative integers. Furthermore, for any ordered set $(I,<)$ we write $(I^*,<_*)$ for the same set with the reversed ordering, i.e.\ for $a,b \in I$ it is $a <_* b :\Leftrightarrow b<a$.
	For a family of ordered sets $(S_i,<_i)_{i \in I}$, we define the \textbf{sum of the orderings} $(\sum_{i \in I} S_i,<)$ as the disjoint union of the $S_i$	with the following ordering: Let $a, b \in I$ and $s_1, s_2 \in \sum_{i \in I} S_i$ with $s_1 \in S_a, s_2 \in S_b$. Then $s_1 < s_2 :\Leftrightarrow [a <_I b \vee (a=b \wedge s_1 <_a s_2)]$. In the case where $I = \{1,\ldots,n\}$ for some $n \in \N$, we also write $\sum_{i \in I} S_i = S_1 + \ldots + S_n$.
	
	Let $(G, +, -, 0, <)$ be an ordered abelian group. For a prime $p \in \N$ we denote the maximal convex $p$-divisible subgroup of $G$ by $G_p$ and the maximal divisible convex subgroup of $G$ by $G_0$. We define the \textbf{absolute value} as $|\cdot|\colon a \mapsto \max\{a, -a\}$.
	We say that two elements $a, b \in G$ are \textbf{archimedean equivalent} and write $a \sim b$ if there exists $n \in \N$ such that $|a| \leq n|b|$ and $|b| \leq n|a|$. The archimedean equivalence class of $a$ in $G$ is given by $[a]_+ := \{ b \in G \mid a \sim b \}$. The set $G/{\sim}  = \{ [a]_+ \mid a \in G \}$ admits a strict ordering given by $[a]_+ < [b]_+ :\Leftrightarrow [|a|>|b| \wedge a \nsim b]$.
	Similar notions are used for an ordered field $(K, +, -, \cdot, 0, 1, <)$ by applying them to its $\Log$-reduct $(K, +, -, 0, <)$.
	
	Let $(I,<)$ be a non-empty ordered set and let $(G_i,<)_{i \in I}$ be a family of ordered abelian groups. We define the \textbf{Hahn product} ${\HProd}_{i \in I}\ G_i$ of $(G_i,<)_{i \in I}$ as the set of all maps $s\colon I \rightarrow \bigsqcup_{i \in I} G_i$ with $s(i) \in G_i$ for any $i \in I$ and with well-ordered \textbf{support} $\supp(s) := \{ i \in I \mid s(i) \in G_i \setminus \{0\} \}$.
	It becomes an ordered abelian group with $(s+r)(i) := s(i) + r(i)$ and the \textbf{lexicographic ordering} $s > 0 :\Leftrightarrow [s \not= 0 \wedge s(\min(\supp(s)))>0]$.
	
	We define the \textbf{Hahn sum} $\coprod_{i \in I} G_i$ as the subgroup of $\HProd_{i \in I} G_i$ consisting of the elements with finite support, i.e.\ 
	$\coprod_{i \in I} G_i=\{s\in \underset{i \in I}{\HProd}\ G_i \mid |\supp(s)|<\infty\}$. Given $n \in \N$, we also write $G_1 \amalg \ldots \amalg G_n$ for $\coprod_{i \in \{1,\ldots,n\}} G_i$. Furthermore, in this case we identify a map $s \in G_1 \amalg \ldots \amalg G_n$ with the $n$-tuple $(s(1),\ldots,s(n))$.
	
	For an ordered field $(k,<)$ and an ordered abelian group $(G,<)$, we denote the (generalised) \textbf{power series field} by $k(\!(G)\!)$. The field $k(\!(G)\!)$ is given by the additive abelian group $\HProd_{g \in G} k$ together with a multiplication defined by
	\[ (s \cdot r)(g) := \sum_{\substack{(h_1,h_2) \in G^2, \\ h_1 + h_2 = g}} s(h_1) \cdot r(h_2). \]
	With these operations and the lexicographic ordering, $(k(\!(G)\!),+,-,\cdot,0,1,<)$ becomes an ordered field (see \cite[page~27~f.]{kuhlmann}).\\
	
	\subsection{Valuations} For a valuation $w$ on a field $K$, we denote its value group by $(w(K^\times), +, -, 0, <)$, its valuation ring $\{ a \in K \mid w(a) \geq 0 \}$ by $\OO_w$ and its valuation ideal $\{ a \in K \mid w(a) > 0 \}$ by $\II_w$. Furthermore, we denote the residue field $\OO_w/\II_w$ by $Kw$. A valuation $w$ on an ordered field $(K,<)$ is called \textbf{convex} if $\OO_w$ is convex in $(K,<)$. We expand an ordered field $(K,+,-,\cdot,<)$ by a fixed valuation $w$ by considering the ordered valued field $(K,+,-,\cdot,<,\OO_w)$, where $\OO_w$ is a unary predicate symbol denoting membership in the valuation ring.
	
	The \textbf{natural valuation} $\vnat$ on an ordered field $(K,<)$ is defined as the finest convex valuation on $(K,<)$. Its valuation ring is given by $\OO_\vnat = \{ a \in K \mid [a]_+ \geq [1]_+ \}$. Note that  $\vnat(a) = \vnat(b) \Leftrightarrow [a]_+ = [b]_+$. Hence for $[a]_+ + [b]_+ := [ab]_+$ we obtain $(K^\times/{\sim},+,-,0,<) \cong (\vnat(K^\times),+,-,0,<)$, see \cite[page 16]{kuhlmann} for details. We often implicitly identify value groups if there is an isomorphism between them. We also note that, by \cite[Lemma 1.1~(i),(iv)]{knebusch}, the convex valuations on $(K,<)$ are exactly the coarsenings of $\vnat$. Let $G := \vnat(K^\times)$. We define $v_p\colon K^\times \longrightarrow G/G_p, x \mapsto \vnat(x) + G_p$.	
	For an ordered abelian group $(G,+,-,0,<)$,	the map $v_G\colon G \rightarrow G/{\sim}, a \mapsto [a]_+$ is the \textbf{natural valuation} on $(G,+,-,0,<)$.
	
	Let $(k,<)$ be an ordered field and let $(G,<)$ be an ordered abelian group. For the power series field $(k(\!(G)\!),<)$, we canonically obtain a valuation $\vmin$ defined by $\vmin(s) := \min(\supp(s))$ for $s \not= 0$. This valuations is convex.
	Since $\vmin$ is always henselian (see \cite[page~27~f.]{kuhlmann}), the field $k(\!(G)\!)$ is real closed if and only if $k$ is real closed and $G$ is divisible (see \cite[Theorem 4.3.7]{engler}).

\section{Definable ranks}\label{section:defRanks}

	\begin{definition}\label{def:Ranks}
		Let $(K,<)$ be an ordered field, $(G,<)$ an ordered abelian group and $(\Gamma,<)$ an ordered set.
		We define the respective \textbf{ranks} of $K, G$ and $\Gamma$ by
		\begin{itemize}
			\item $\rk_K := \{ \OO_v \mid \OO_v \textrm{ convex valuation ring on } K \textrm{ with } \OO_\vnat \subseteq \OO_v \subsetneq K \}$,
			
			\item $\rk_G := \{ H \mid H \textrm{ convex subgroup of } G \textrm{ with } \{0\} \subseteq H \subsetneq G \}$,
			
			\item $\rk_\Gamma := \{ \Delta \mid \Delta \textrm{ final segment of } \Gamma \textrm{ with } \emptyset \subseteq \Delta \subsetneq \Gamma \}$.
		\end{itemize}
	\end{definition}

	Note that $\rk_K, \rk_G$ and $\rk_\Gamma$ are each ordered by $\subsetneq$.
	
	\begin{remark}
		\begin{itemize}
			\item[(a)] The set $\rk_K$ consists exactly of the valuation rings of non-trivial coarsenings of $\vnat$.
			
			\item[(b)] In the literature one can find different notions of the rank. The notion given here is similar to the one found in \cite[page~26--28]{engler}. In contrast, the definition of rank from \cite[page~50]{kuhlmann} excludes the minimal elements $\OO_\vnat, \{0\}$ and $\emptyset$ respectively, but includes the entire set as maximal element unless the minimal element was already the entire set. The exclusion of exactly one minimum and maximum is done to obtain rank $1$ exactly for archimedean ordered groups and rank 0 for archimedean ordered fields. Usually only the order type of the sets above is considered when one talks about the rank. We however need to really look at the exact sets for reasons that will become apparent shortly.
		\end{itemize}
	\end{remark}

	The following establishes the relation between the rank of an ordered field, that of its value group under the natural valuation and the rank of the value set of the natural valuation on the former group (see \cite[Lemma 3.4,~3.5]{kuhlmann} for a proof).

	\begin{fact}\label{fct:rankIsomorphy}
		Let $(K,<)$ be an ordered field. Denote $G = \vnat(K^\times)$ and $\Gamma = v_G(G\setminus\{0\})$. Then
		\[ \begin{array}{ccc}
			\rk_K & \longrightarrow & \rk_G\\
			\OO_v & \longmapsto & \vnat(\OO_v \setminus \II_v)
		\end{array} \]
		is an isomorphism from $(\rk_K, \subsetneq)$ to $(\rk_G, \subsetneq)$. Analogously
		\[ \begin{array}{ccc}
			\rk_G & \longrightarrow & \rk_\Gamma\\
			H & \longmapsto & v_G(H\setminus\{0\})
		\end{array} \]
		is an isomorphism from $(\rk_G, \subsetneq)$ to $(\rk_\Gamma, \subsetneq)$.
	\end{fact}	

	\begin{notation}
		We denote the above maps by $\Phi_K\colon \OO_v \mapsto \vnat(\OO_v \setminus \II_v)$ and $\Phi_G\colon H \mapsto v_G(H\setminus\{0\})$ respectively.
	\end{notation}

	\begin{definition}\label{def:definableRank}
		Let $(K,<)$ be an ordered field, $(G,<)$ an ordered abelian group and $(\Gamma,<)$ an ordered set.
		We define the respective \textbf{definable ranks} of $K, G$ and $\Gamma$ as
		\begin{itemize}
			\item $\drk_K := \{ \mathcal{O} \in \rk_K \mid \mathcal{O} \textrm{ is } \mathcal{L}_{\mathrm{or}}\textrm{-definable in } K \}$,
		
			\item $\drk_G := \{ H \in \rk_G \mid H \textrm{ is } \mathcal{L}_{\mathrm{og}}\textrm{-definable in } G  \}$,
		
			\item $\drk_\Gamma :=  \{ \Delta \in \rk_\Gamma \mid \Delta \textrm{ is } \mathcal{L}_{<}\textrm{-definable in } \Gamma  \}$.
		\end{itemize}
	\end{definition}

	While \autoref{fct:rankIsomorphy} establishes the pairwise equivalence of $\rk_K, \rk_G$ and $\rk_\Gamma$ as ordered sets, we give an example, which shows that their definable counterparts are not necessarily pairwise isomorphic.
	
	\begin{example}\label{xmpl:ArchG}
		Let $K$ be an ordered field such that $G = \vnat(K^\times)$ is ar\-chi\-me\-de\-an and denote $\Gamma = v_G(G \setminus \{0\})$. Then $\rk_K = \{\OO_\vnat\}, \rk_G = \{ \{0\} \}$ and $\rk_\Gamma = \{\emptyset\}$. It follows immediately that $\drk_G = \{ \{0\} \}$ and $\drk_\Gamma = \{\emptyset\}$.
		
		In order for all three definable ranks to coincide, $\vnat$ must now be definable in $K$.
		\begin{itemize}
			\item [(i)] If $G$ is not divisible or $K\vnat$ is not real closed, then by \cite[Corollary 3.2]{DittmannJahnkeKrappKuhlmann23} it follows that $\OO_\vnat$ is definable and all three definable ranks coincide.
			
			\item[(ii)] Let $K$ be the field of formal Puiseux series $K = \bigcup_{n \in \NN} \R(\!(t^{1/n})\!)$, ordered by $0 < t < \Q^{>0}$. Then $K$ is real closed and $\vnat$ is not definable due to o-minimality. Hence $\drk_K = \emptyset$ and the definable rank on the field-level does not coincide with the definable ranks on the value group- and value set-level.
		\end{itemize}
	\end{example}

	\begin{example}\label{xmpl:allRanksDifferent}
		Let $K := \R(\!(\Q \amalg \Q)\!)$. Then $K$ is real closed, hence $\drk_K = \emptyset$ as in \autoref{xmpl:ArchG}~(ii). Furthermore $G := \vnat(K^\times) = \Q \amalg \Q$ is divisible, by o-minimality it follows $\drk_G = \{ \{0\} \}$. Lastly, $\Gamma := v_G(G \setminus \{0\})$ has two elements, write $\Gamma = \{ 1, 2 \}$. Then $\drk_{\Gamma} = \{ \emptyset, \{2\} \}$. We see that in this example all definable ranks are non-isomorphic as ordered sets.
	\end{example}

	\begin{example}\label{xmpl:drkNonIsomorphic}
		Consider the field $k := \R(t^{\frac{1}{n}} \mid n \in \N) \subsetneq \R(\!(\Q)\!)$. We can consider $k$ as an ordered field with the restriction of the unique ordering of $\R(\!(\Q)\!)$. Let $\vnat'$ be the natural valuation on $(k,<)$. As shown in \cite[Example 3.7]{DittmannJahnkeKrappKuhlmann23}, the natural valuation is $\Lor$-definable by a formula $\varphi(x)$, $k\vnat' = \R$ and $\vnat'(k^\times) = \Q$. Since $\vnat$ is the only non-trivial convex valuation on $(k,<)$, it follows that the $\Lor$-theory of $(k,<)$ contains sentences expressing that every non-trivial definable valuation ring is equal to the one defined by $\varphi(x)$. We now consider an elementary extension $K \succ k$ such that the partial type \[ \Sigma(x) := \{ x > n \mid n \in \Z \} \cup \{ \varphi(x) \} \] is realized in $(K,<)$. Then $\varphi(K)$ is the only definable non-trivial convex valuation ring of $(K,<)$, denote the corresponding valuation by $w$. Now $w$ is not the natural valuation on $(K,<)$ since the residue of a realization of $\Sigma(x)$ is greater than all integers, making $Kw$ non-archimedean.
		
		On the other hand, the $\Lor$-theory of $(k,<)$ also contains sentences expressing that the valuation ring defined by $\varphi(x)$ has a divisible value group and a real closed residue field. The natural valuation on $(K,<)$ is the composition of $w$ and the natural valuation on $(Kw,<)$. Since $Kw$ must be real closed, it follows from \cite[Proposition 4.3.7]{engler} it has a divisible value group. It now follows that the composition of two valuations with divisible value groups must itself have a divisible value group. Hence $G := \vnat(K^\times)$ is divisible and its only definable proper convex subgroup is the trivial one.
		
		Together we obtain $\drk_K = \{\OO_w\}$ and $\drk_G = \{ \{0\} \}$, but $\Phi_K(\OO_w) \not= \{0\}$. So while the definable ranks of $K$ and $G$ are isomorphic as ordered sets, they are not in correspondence under $\Phi_K$ and in fact $\Phi_K(\drk_K) \cap \drk_G = \emptyset$.
	\end{example}
	
\section{Definability of convex subgroups}\label{section:defCSubgroups}
	We now want to discuss when a given convex subgroup $H$ of an ordered abelian group $(G,+,<)$ is $\Log$-definable. To do so, we heavily rely on Schmitt's work on the model theory of ordered abelian groups (see \cite{Schmitt1982}). We ultimately obtain a full characterization of definable convex subgroups, which is effectively \cite[Corollary 4.2]{delon}; we are going to show that the converse implication of the one in the corollary holds as well. We then conclude the section by illustrating that in many cases it suffices to check $p$-divisibility of certain group quotients for primes $p \in \N$, allowing us to deduce definability results on purely algebraic properties.
	
	At first we need to elaborate on the necessary terminology and definitions to actually employ the results from \cite{Schmitt1982}. We start with the following property, which was first coined by Robinson and Zakon in \cite[Definition 3.3.~(v)]{robinson} and later generalized to only consider a fixed integer $n$.

	\begin{definition}\label{def:nRegular}
		Let $n \in \N$ with $n \geq 2$. An ordered abelian group $G$ is called $n$\emph{-regular} if for every interval $I := [a,b] \subsetneq G$ with $|I| \geq n$ there exists $t \in G$ such that $nt \in I$.
	\end{definition}

	We note that this property is certainly expressible as a first order sentence in the language $\Lor$. We also cite the following alternative characterization of $n$-regular groups.

	\begin{fact}\label{fact:nRegChar}
		\textrm{\cite[Lemma 1.14.]{Schmitt1982}} A non-trivial ordered abelian group $G$ is $n$-regular if and only if for every non-trivial convex subgroup $H \subseteq G$ the quotient $G/H$ is $n$-divisible.
	\end{fact}

	As pointed out in \cite[page~236]{robinson}, the regular ordered abelian groups are exactly those that are elementary equivalent to archimedean ordered abelian groups. Since the only proper convex subgroup of an archimedean ordered abelian group is the trivial one, there is no hope to define any non-trivial proper convex subgroup of a regular ordered abelian group. Now we observe that every formula can only consider $n$-regularity for finitely many $n \in \N$ (which is equivalent to only considering $\ell$-regularity where $\ell$ is their least common multiple). While in the classical setting one now deconstructs an ordered abelian group into a collection of archimedean ordered abelian groups, in the context of model theory one should now instead break them down into regular ordered abelian groups. Since regularity itself is not a single first order property, Schmitt instead uses $n$-regularity for arbitrary $n \in \N$ and considers them separately.

	\begin{definition}\label{def:SchmittSubgroups}
		Let $G$ be an ordered abelian group, let $g \in G\setminus \{0\}$ and let $n \in \N, n \geq 2$. Further denote $S(G) := \{H \subseteq G \mid H \textrm{ convex subgroup} \}$. We define
		\[ A(g) := \bigcup \{ H \in S(G) \mid g \not\in H \}, \]
		\[ B(g) := \bigcap \{ H \in S(G) \mid g \in H \}, \]
		\[ C(g) := B(g)/A(g), \]
		\[ A_n(g) := \bigcap \{ H \in S(G) \mid H \subsetneq B(g), B(g)/H \textrm{ is }n\textrm{-regular} \}, \]
		\[ B_n(g) := \bigcup \{ H \in S(G) \mid A(g) \subsetneq H, H/A(g) \textrm{ is }n\textrm{-regular} \}\textrm{ and} \]
		\[ C_n(g) := B_n(g)/A_n(g). \]
		We furthermore define 
		\[ F_n(g) := \bigcup \{ H \in S(G) \mid H\cap(g+nG) = \emptyset \} \] for $g \not\in nG$ and $F_n(g) := \emptyset$ for $g \in nG$. Lastly, we can now define
		\[ E_n(g) := \{ h \in G \mid F_n(h) \subseteq F_n(g) \}, \]
		\[ E_n^*(g) := \{ h \in G \mid F_n(h) \subsetneq F_n(g) \} \textrm{ and} \]
		\[ F_n^*(g) := E_n(g)/E_n^*(g). \]
	\end{definition}


	\begin{definition}\label{def:nSpines}
		\begin{itemize}
			\item[(1)] The \emph{language of spines} $\Lsp$ consists of a binary predicate symbol $\geq$ and the  unary predicate symbols $A$, $F$, $\mathrm{Dk}$ and $\alpha(p,k,m)$ for all primes $p \in \N$, $k \in \N$ and $m \in \omega$.
			
			\item[(2)] Let $G$ be an ordered abelian group and $n \in \N, n \geq 2$. The $n$\emph{-spine of} $G$, denoted as $\SPn(G)$, is the $\Lsp$-structure with the domain \[ \{A_n(g) \mid g \in G \setminus \{0\} \} \cup \{ F_n(g) \mid g \in G \setminus nG \} \] and the following interpretation of the symbols from $\Lsp$:
			\begin{itemize}
				\item[(i)] $[\SPn(G) \models C \geq D] \Leftrightarrow C \subseteq D$,
				\item[(ii)] $[\SPn(G) \models A(C)] \Leftrightarrow \exists g \in G \setminus \{0\}\colon C = A_n(g)$,
				\item[(iii)] $[\SPn(G) \models F(C)] \Leftrightarrow \exists g \in G \setminus nG\colon C = F_n(g)$,
				\item[(iv)] $[\SPn(G) \models \mathrm{Dk}(C)] \Leftrightarrow G/C \textrm{ is discrete}$,
				\item[(v)] $[\SPn(G) \models \alpha(p,k,m)(C)] :\Leftrightarrow
				\begin{array}{c}
					\exists g \in G \setminus nG\colon C = F_n(g) \\ \wedge \dim(p^kF_n^*(g)/p^{k+1}F_n^*(g)) \geq m.
				\end{array}$
			\end{itemize}
		\end{itemize}
	\end{definition}

	The main result of Schmitt is a quantifier elimination, essentially up to a formula over one of the $n$-spines. The result is quite complex and lengthy to formulate, for our purposes it suffices that the converse direction also holds. Every $\Lsp$-formula can be translated into an $\Log$-formula, in a way that is made precise in the following.

	\begin{fact}\label{fact:reverseGroupQE}
		\textrm{\cite[Lemma 3.1.]{Schmitt1982}} Let $n \geq 2$, let $\varphi(y_1,\ldots,y_k,z_1,\ldots,z_\ell)$ be an $\Lsp$-formula and let $t_i(x_1,\ldots,x_m)$ and $s_j(x_1,\ldots,x_m)$ be $\Log$-terms, $1 \leq i \leq k$ and $1 \leq j \leq \ell$. Then there is an $\Log$-formula $\bar{\varphi}(x_1,\ldots,x_m)$ such that for all ordered abelian groups $G$ and all $g_1,\ldots,g_m \in G$ it holds that\smallskip\\
		$G \models \bar{\varphi}(g_1,\ldots,g_m)$ if and only if
		\[ \SPn(G) \models \varphi(C_1,\ldots,C_k,D_1,\ldots,D_\ell) \]
		where $C_i := A_n(t_i(g_1,\ldots,g_m))$ and $D_j := F_n(s_j(g_1,\ldots,g_m))$.
	\end{fact}

	To properly apply Schmitt's theory to our purpose, we furthermore use a variety of results from his work. The results which are needed are collected below for the readers convenience.

	\begin{fact}\label{fact:SchmittCollection}
		Let $n \in \N, n \geq 2$, let $G$ be an ordered abelian group and $g,h \in G$. It holds
		\begin{itemize}
			\item[(1)] \cite[Lemma 2.1.~(3)]{Schmitt1982} $A_n(g) \subseteq A(g) \subsetneq B(g) \subseteq B_n(g)$,
			
			\item[(2)] \cite[Lemma 2.2.~(1)]{Schmitt1982} If $A_n(g) \subsetneq A_n(h)$, then $B_n(g) \subseteq A_n(h)$ and
			
			\item[(3)] \cite[Lemma 2.9.~(3)]{Schmitt1982} $F_n(g) = \bigcap\{ A_n(g+nd) \mid d \in G \}$.
		\end{itemize}
	\end{fact}

	We now have all the required ingredients in place to fully characterize when a convex subgroup $H \subseteq G$ is $\Log$-definable in terms of the $n$-spines of $G$. Note that one direction of the characterization can be found in \cite{delon}, the backwards direction however is never proved. We fill in the missing direction and furthermore make sure that nothing is lost through the multitude of notational differences across the literature.

	\begin{theorem}\label{theorem:CharacterizationOfDefinableConvexSubgroups}
		Let $(G,+,0,<)$ be an ordered abelian group and $H \subsetneq G$ a convex subgroup. $H$ is $\Log$-definable if and only if there is an $n \in \N, n \geq 2$ and an $\Lsp$-definable initial segment $\Delta \subseteq \SPn(G)$ such that \[ H = \bigcap_{C \in \Delta} C. \]
	\end{theorem}
	\begin{proof}
		`$\Rightarrow$': Let $H$ be $\Log$-definable, then by \cite[Corollary 4.2.]{delon} there is an $n \geq 2$ and an $\Lsp$-definable final segment $\Delta'$ (note that in \cite{delon} spines are considered with reverse ordering) such that $g \in H$ if and only if $A_n(g) \in \Delta'$ for a $g \in G$. Consider $\Delta := \SPn(G) \setminus \Delta'$. By \cite[Theorem 4.1.]{delon} for some $n$ it is $H = \bigcap_{g \not\in H} A_n(g)$ and following the proof of  \cite[Corollary 4.2.]{delon}, the $n$ is the same. Then $H = \bigcap_{g \not\in H} A_n(g) = \bigcap_{C \in \Delta, \SPn(G) \models A(C)} C = \bigcap_{C \in \Delta} C$, where the last equation follows since for all $C \in \Delta$ with $\SPn(G) \not\models A(C)$ it follows $\SPn(G) \models F(C)$ and by \autoref{fact:SchmittCollection}~(3) $C$ is the intersection of $C' \leq C$ with $\SPn(G) \models A(C')$.\\
	
		`$\Leftarrow$': Fix $n \geq 2$, an initial segment $\Delta \subseteq \SPn(G)$ and an $\Lsp$-formula $\varphi(y)$ such that $\varphi(\SPn(G)) = \Delta$. Let $H := \bigcap_{C \in \Delta} C$, it remains to show that $H$ is $\Log$-definable in $G$. Consider $\psi(y) := A(y) \wedge \neg\varphi(y)$ and the $\Log$-term $t(x) = x$. By \autoref{fact:reverseGroupQE} there is an $\Log$-formula $\bar{\psi}(x)$ such that $G \models \bar{\psi}(g) \Leftrightarrow \SPn(G) \models \psi(A_n(g))$ for all $g \in G$. It suffices to show that $H = \bar{\psi}(G)$.
		\begin{itemize}
			\item[(1)] Let $g \in H$, then $A_n(g) \subsetneq H$ since $g \not\in A_n(g)$ and thus $A_n(g) \not\in \Delta$. Therefore $\SPn(G) \models \psi(A_n(g))$ and $g \in \bar{\psi}(G)$.
			
			\item[(2)] Let $g \in G \setminus H$. Then there is $C_0 \in \Delta$ with $g \not\in C_0$. By \autoref{fact:SchmittCollection}~(3) we can without loss of generality assume that $\SPn(G) \models A(C_0)$. Now since $g \not\in C_0$, it follows that $C_0 \subseteq A(g) \subsetneq B(g) \subseteq B_n(g)$ by \autoref{fact:SchmittCollection}~(1). Furthermore $C_0 \subseteq A_n(g)$, because $C_0 = A_n(h)$ for some $h \in G$ and otherwise we had $A_n(g) \subsetneq A_n(h) \subsetneq B_n(g)$, contradicting \autoref{fact:SchmittCollection}~(2). Therefore $A_n(g) \leq C_0$ and thus $A_n(g) \in \Delta$ since $\Delta$ is an initial segment. This yields $\SPn(G) \not\models \psi(A_n(g))$ and $g \not\in \bar{\psi}(G)$.
		\end{itemize}
		We have therefore shown  that $g \in H$ if and only if $g \in \bar{\psi}(G)$ for all $g \in G$, completing the proof.
	\end{proof}

	With this we have established a conclusive characterization of definable convex subgroups $H \subseteq G$. However, we did so by utilizing a different model theoretic structure. In a multitude of cases it turns out that purely order theoretic and algebraic properties of $G, H$ and certain quotients of convex subgroups are sufficient to decide whether $H$ is $\Log$-definable. We now give several sufficient conditions and deduce them from our characterization.

	\begin{corollary}\label{coro:dfblConvSubgroups}
		Let $G$ be an ordered abelian group and $H \subsetneq G$ a convex subgroup.
		\begin{itemize}
			\item[(1)] If $\Phi_G(G) \setminus \Phi_G(H)$ has a maximum and there is an $n \geq 2$ such that $G/H$ has no non-trivial $n$-divisible convex subgroup, then $H$ is $\Log$-definable.
			
			\item[(2)] If for some $n \geq 2$ the quotient $G/H$ has a non-trivial $n$-divisible convex subgroup, $\Phi_G(H)$ does not have a minimum and for every convex subgroup $H' \subsetneq H$ the quotient $H/H'$ is not $n$-divisible, then $H$ is $\Log$-definable.
			
			\item[(3)] If for some $n \geq 2$ the quotient $G/H$ has no non-tivial convex $n$-divisible subgroup and $\Phi_G(H)$ has a minimum, then $H$ is $\Log$-definable.
			
			\item[(4)] If for some $n \geq 2$ the quotient $G/H$ has no non-trivial convex $n$-divisible subgroup and for some convex subgroup $H' \subsetneq H$ the quotient $H/H'$ is $n$-divisible, then $H$ is $\Log$-definable.
		\end{itemize}
	\end{corollary}
	\begin{proof}
		We show in all cases that $H = \bigcap_{C \in \Delta} C$ for $\Delta \subseteq \SPn(G)$ a definable initial segment.
		\begin{itemize}
			\item[(1)] Consider $g \in G$ with $v_G(g) = \max(\Phi_G(G) \setminus \Phi_G(H))$. Then $H = A(g)$. Since $G/H$, it follows that $H \subseteq A_n(g)$, so $H = A_n(G)$ by \autoref{fact:SchmittCollection}~(1). Then $H = \bigcap_{C \leq A_n(g)} C$ and $\Delta := \{ C \in \SPn(G) \mid C \leq A_n(g) \}$ is a definable initial segment of $\SPn(G)$.
			
			\item[(2)] Let $G' \subseteq G$ be a convex subgroup with $H \subsetneq G'$ such that $G'/H$ is $n$-divisible. Choose $g \in G' \setminus H$, then $A_n(g) \subseteq H$ as $B(g)/H$ is a convex subgroup of $G'/H$ and thus $n$-divisible and in particular $n$-regular.\newline Assume $A_n(g) \subsetneq H$, then $H/A_n(g)$ is a convex subgroup of $B(g)/A_n(G)$ and thus $n$-regular. Therefore for every convex $H' \subsetneq G$ with $A_n(g) \subsetneq H' \subsetneq H$ it follows that $(H/A_n(g))/(H'/A_n(g)) \cong H/H'$ is divisible, which contradicts the conditions. Hence there is no such $H'$, but then $H = B(h)$ for any $h \in H \setminus A_n(g)$ and $v_G(h)$ is minimal in $\Phi_G(H)$. Contradiction to the conditions. Therefore $H = A_n(G)$ and $H = \bigcap_{C \leq A_n(g)} C$ as in (1).
			
			\item[(3)] Consider $h \in H$ with $v_G(h) = \min(\Phi_G(H))$. Then $A(h) \subsetneq H = B(h)$. Now assume $H \subsetneq B_n(h)$. Then $B_n(h)/A(h)$ is $n$-regular and $H/A(h) \subsetneq B_n(h)/A(h)$ is a non-trivial convex subgroup.\newline
			Thus $(B_n(h)/A(h))/(H/A(h)) \cong B_n(h)/H$ is $n$-divisible. Contradiction to the conditions. It follows that $H = B_n(h)$ and for all $g \in G \setminus H$ we obtain $H \subseteq A_n(g)$. On the other hand $\bigcap_{g \in G \setminus H} A_n(g) \subseteq H$, so $H = \bigcap_{C \in \Delta} C$ for $\Delta := \{ C \in \SPn(G) \mid C < A_n(h) \}$, which is a definable initial segment of $\SPn(G)$.
			
			\item[(4)] Let $H' \subsetneq H$ be a convex subgroup such that $H/H'$ is $n$-divisible. Choose $h \in H \setminus H'$, then $H' \subseteq A(h)$ and $H/A(h) \cong (H/H')/(A(h)/H')$ is $n$-divisible. Therefore $H \subseteq B_n(h)$. Assume $H \subsetneq B_n(h)$, then $B_n(h)/H$ is a non-trivial $n$-divisible convex subgroup of $G/H$, which contradicts the conditions. This implies $H = \bigcap_{C \in \Delta} C$ for $\Delta := \{ C \in \SPn(G) \mid C < A_n(h) \}$ as in (3).
		\end{itemize}
		The $\Log$-definability of $H$ now follows from \autoref{theorem:CharacterizationOfDefinableConvexSubgroups} in all cases.
	\end{proof}

	\begin{corollary}\label{coro:ndfblConvSubgroups}
		Let $G$ be an ordered abelian group and $H \subsetneq G$ a convex subgroup. If for every $n \geq 2$ the quotient $G/H$ has a non-trivial $n$-divisible convex subgroup and
		\begin{itemize}
			\item[(a)] $\Phi_G(H)$ has a minimum or
			
			\item[(b)] there is a convex subgroup $H' \subsetneq H$ such that $H/H'$ is $n$-divisible,
		\end{itemize}
		then $H$ is not $\Log$-definable.
	\end{corollary}
	\begin{proof}
		Assume for contradiction that $H \subsetneq G$ as above was $\Lor$-definable. Then by \autoref{theorem:CharacterizationOfDefinableConvexSubgroups} there is an $n \geq 2$ and an $\Lsp$-definable initial segment $\Delta$ of $\SPn(G)$ such that $H = \bigcap_{C \in \Delta} C$. By the assumptions, there is a convex subgroup $G' \subseteq G$ such that $H \subsetneq G'$ and $G'/H$ is $n$-divisible. Choose $h \in G' \setminus H$, then $H \subsetneq B(h) \subseteq G'$ and since $B(h)/H$ is $n$-divisible, as it is a convex subgroup of $G'/H$, it follows $A_n(h) \subseteq H$.
		
		Case (a): If $\Phi_G(H)$ has a minimum $\gamma$, choose $h' \in H$ with $v_G(h') = \gamma$. Note that $H = B(h')$. Now $B(h)/A(h')$ is $n$- regular because every non-trivial convex subgroup of $B(h)/A(h')$ contains $h' + A(h')$ and thus has the form $C/A(h')$ where $B(h') \subseteq C \subseteq B(h)$. But then $(B(h)/A(h'))/(C/A(h')) \cong B(h)/C \cong (B(h)/H)/(C/H)$. This is $n$-divisible, as it is a quotient of the $n$-divisible group $B(h)/H$. Therefore $A_n(h) \subseteq A(h')$ and $A_n(h) \subsetneq H$.
		
		Case (b): If for some convex subgroup $H' \subsetneq H$ the quotient $H/H'$ is $n$-divisible, then $B(h)/H'$ is $n$-divisible. Therefore $A_n(h) \subseteq H'$ and $A_n(h) \subsetneq H$.
		
		Since we assumed one of the cases, we obtain that $A_n(h) \subsetneq H$. Hence for all $C \in \Delta$ it certainly is $C < A_n(h)$. However, for all $C \in \Delta$ with $\SPn(G) \models A(C)$ it now follows $C \subseteq B_n(h) \supseteq B(h)$ from \autoref{fact:SchmittCollection}~(2). In particular $h \in C$ for such $C$. Furthermore, for every $C' \in \Delta$ with $\SPn(G) \models F(C')$, it follows from \autoref{fact:SchmittCollection}~(3) that $C'$ is the intersection of $C \in \Delta$ with $\SPn(G) \models A(C)$. By the prior this also yields $h \in C'$. But then $h \in \bigcap_{C \in \Delta} C$ and thus $\bigcap_{C \in \Delta} C \not= H$. Contradiction. Hence $H$ was not $\Lor$-definable.
	\end{proof}

	We now conclude this section by giving a sufficient condition on an ordered abelian group $G$ for every definable final segment of the value set $\Gamma := v_G(G \setminus \{0\})$ to correspond to a definable convex subgroup $H \subseteq G$. Note that this condition is not necessary.

	\begin{proposition}\label{prop:suffConditionFordkrGContainsdrkGamma}
		Let $G$ be an ordered abelian group and $\Gamma = v_G(G \setminus \{0\})$. If for some $n \in \N, n \geq 2$ and all $g \in G \setminus \{0\}$ the groups $C(g)$ are not $n$-divisible, then $\drk_{\Gamma} \subseteq \Phi_G(\drk_G)$.
	\end{proposition}
	\begin{proof}
		Let $g \in G \setminus \{0\}$, $n \geq 2$ as above. Then $A_n(g) = A(g)$, otherwise $A_n(g) \subsetneq A(g)$ and $(B(g)/A_n(g))/(A(g)/(A_n(g))) \cong B(g)/A(g) \cong C(g)$ is not $n$-divisible and thus $B(g)/A_n(g)$ not $n$-regular. Contradiction.
		
		Now $A(g)$ only depends on $v_G(g)$ for every $g \in G \setminus \{0\}$ and for $g,h \in G \setminus \{0\}$ with $v_G(g) > v_G(h)$ we obtain $A(g) \subsetneq A(h)$. Hence the map
		\[\begin{array}{cccc}
			\Theta\colon & \{ C \in \SPn(G) \mid \SPn(G) \models A(C) \} & \longrightarrow & \Gamma\\
			& A_n(g) & \longmapsto & v_G(g)
		\end{array}\]
		is an order-reversing bijection. Thus if $\Delta \subseteq \Gamma$ is an $\L_{\mathrm{<}}$-definable final segment with defining formula $\varphi(x)$, then $\Theta^{-1}(\Delta)$ is an $\Lsp$-definable initial segment of $\{ C \in \SPn(G) \mid \SPn(G) \models A(C) \}$ by reversing every ordering in $\varphi(x)$, replacing every parameter $\gamma \in \Gamma$ with $\Theta^{-1}(\gamma)$ and conjoining with $A(x)$ and $A(y)$ for every $y$ appearing in a quantifier. Denote the $\Lsp$-formula obtained in this way as $\tilde{\varphi}(x)$ and consider $\psi(x) := \exists y\colon \tilde{\varphi}(y) \wedge x \leq y$. With $H_\Delta := \bigcap_{C \in \psi(\SPn(G))} C$ we now obtain that $\Phi_G^{-1}(\Delta) = H_\Delta \in \drk_G$ for all $\Delta \in \drk_{\Gamma}$.
	\end{proof}

\section{Comparing group- and field-level}\label{section:FieldAndGroup}
	In this section we apply results from \cite{DittmannJahnkeKrappKuhlmann23} and \cite{KrappKuhlmannLink23} to compare the definable rank of an ordered field $(K,<)$ and that of the value group of its natural valuation $G := \vnat(K^\times)$. This leads to our main result \autoref{theorem:MainTheorem}, which conclusively compares the two definable ranks for all ordered fields with a henselian natural valuation. The crucial ingredient is the following result on ordered abelian groups, which allows us to apply \cite[Theorem 3.1]{DittmannJahnkeKrappKuhlmann23} to several definable convex subgroups of $G$.

	\begin{proposition}\label{prop:convexSubGKeyLemma}
		Let $G$ be an ordered abelian group. Consider a convex subgroup $H \subseteq G$ such that $G_0 \subsetneq H$. Then there are convex subgroups $H_a, H_b \subseteq G$ such that $G_0 \subseteq H_b \subsetneq H_a \subseteq H$ and $H_a/H_b$ is discrete or not closed in its divisible hull.
	\end{proposition}
	\begin{proof}
		First note that $H/G_0$ does not have any non-trivial convex divisible subgroup by the definition of $G_0$. In particular, there is $n \in \N$ and $g \in H$ such that $g + G_0$ is not $n$-divisible. Note that this implies $G_0 \cap g + nG = \emptyset$ and therefore $G_0 \subseteq F_n(g) \subsetneq B(g)$. Now let $H_a := B(g)$ and $H_b := F_n(g)$. If $H_a/H_b$ is discrete we are done, so assume $H_a/H_b$ to be dense.
	
		Case 1: $H_a/H_b$ has a smallest non-trivial convex subgroup $H'/H_b$.
	
		Then $H'$ is the smallest convex subgroup og $G$ such that $H_b \subsetneq H'$. Note that it follows $H' \cap g + nG \not= \emptyset$ since $H_b = F_n(g) \subsetneq H'$, hence we find $h \in H' \setminus H_b$ such that $h + nG = g + nG$ and $F_n(g) = F_n(h)$ by \autoref{fact:SchmittCollection}~(3). Furthermore $H_b \subsetneq B(h) \subseteq H'$ and since $H'$ was minimal we obtain $B(h) = H'$ and $A(h) = H_b$. Therefore $H'/H_b$ is archimedean, which implies that $H'/H_b$ is $n$-regular. On the other hand $H'/H_b$ is not $n$-divisible since $h$ can not be divided by $n$. Because $H'/H_b \subseteq H_a/H_b$ is a convex subgroup that is $n$-regular and not $n$-divisible, it follows from \cite[Proposition 3.3]{KrappKuhlmannLink23} that $H_a/H_b$ is not closed in its divisible hull.
	
		Case 2: $H_a/H_b$ has no smallest non-trivial convex subgroup.
	
		We show that $(g + H_b)/n$ is a limit point of $H_a/H_b$ in the divisible hull, i.e.\ every open interval $I \subseteq (H_a/H_b)^\mathrm{div}$ that contains $(g + H_b)/n$ already contains a point from $H_a/H_b$. It suffices to show that $nI \cap H_a/H_b$ contains an $h \in n(H_a/H_b)$, since then $h/n \in I \cap H_a/H_b$.
	
		Assume for contradiction there was $\delta \in H_a \setminus H_b$ with $\delta+H_b > 0$ such that $((g-\delta)+H_b, (g+\delta)+H_b) \cap n(H_a/H_b) = \emptyset$.  Now $A(\delta) \subsetneq H_a$. Furthermore $H_b \subsetneq A(\delta)$, since otherwise $B(\delta)/H_b ) C(\delta)$ would be an archimedean non-trivial convex subgroup of $H_a/H_b$. Then it would be the smallest non-trivial convex subgroup of $H_a/H_b$, but the case distinction asserts that this does not exist. Therefore $A(\delta) \supsetneq F_n(g)$, i.e.\ $A(\delta) \cap g + nG \not= \emptyset$. Thus there is $h \in A(\delta)$ such that $g+h$ is $n$-divisible in $G$ and therefore in $H_a$. But $g+A(\delta) \subsetneq (g-\delta, g+\delta)$, hence $(g+h) + H_b \in ((g-\delta)+H_b, (g+\delta)+H_b)$ and $n$-divisible. Contradiction.
	
		Therefore if $H_a/H_b$ is densely ordered, then it is not closed in its divisible hull.
	\end{proof}

	For all convex subgroups $H \in \drk_G$ that are greater than $G_0$, we now obtain that \cite[Theorem 3.1.~(i),(ii)]{DittmannJahnkeKrappKuhlmann23} is applicable to $H$ or an $H' \subseteq H$. This now yields for all such $H$ that $\{\Phi_K\}^{-1}(H)$ is a coarsening of a definable valuation ring.

	\begin{lemma}\label{lem:CoarseningOfDefinableWithDfblSubgp}
		Let $(K,<)$ be an ordered field and let $w$ be an $\L_{\mathrm{or}}$-definable convex valuation on $K$. Let further $H \subseteq w(K^\times)$ be an $\L_{\mathrm{og}}$-definable convex subgroup. Then the coarsening $v$ of $w$ with $\OO_v = \{a \in K \mid w(a) > 0 \vee w(a) \in H \}$ is an $\L_{\mathrm{or}}$-definable valuation on $K$.
	\end{lemma}
	\begin{proof}
		Since $w$ is $\L_{\mathrm{or}}$-definable, its value group is interpretable in the field $(K,+,-,\cdot,<)$. Because $H$ was $\L_{\mathrm{og}}$-definable in the value group of $w$, we find an $\L_{\mathrm{or}}$-formula $\varphi(x)$ such that for $a \in K$, $K \models \varphi(a)$ if and only if $w(a) \in H$. Now $\OO_v = \varphi(K) \cup \OO_w$ is the union of two definable sets and thus definable.
	\end{proof}	

	We can now combine these results to obtain the following.

	\begin{theorem}\label{theorem:mostDfblSubGpCorrespDfblVal}
		Let $(K,<)$ be an ordered field and let $G := \vnat(K^\times)$. Let $H \subseteq G$ be an $\Log$-definable convex subgroup of $G$ such that $G_0 \subsetneq H$. Then the coarsening $\Phi_K^{-1}(H)$ of $\vnat$ is $\Lor$-definable in $K$.
	\end{theorem}
	\begin{proof}
		By \autoref{prop:convexSubGKeyLemma} there are convex subgroups $H_a,H_b \subseteq H$ such that $H_a/H_b$ is discrete or not closed in its divisible hull. In particular, $G/H_b$ is discrete or not closed in its divisible hull and for $v := \Phi_K^{-1}(H_b)$ we have $v(K^\times) = G/H_b$. Then by \cite[Theorem 3.1]{DittmannJahnkeKrappKuhlmann23} it follows that $v$ is an $\Lor$-definable convex valuation. Now $\Phi_K^{-1}(H)$ is a coarsening of $v$ and by \cite[Proposition 3-1]{DelonLucas1989} it follows that $H/H_b$ is $\Log$-definable in $G/H_b = v(K^\times)$. Now \autoref{lem:CoarseningOfDefinableWithDfblSubgp} implies the $\Lor$-definability of $\Phi_K^{-1}(H)$.
	\end{proof}

	In particular, from the above we immediately obtain that for any ordered field $(K,<)$ and $G := \vnat(K^\times)$, only few convex subgroups of $G$ may be $\Log$-definable without corresponding to an $\Lor$-definable convex valuation on $K$. We in fact can deduce the following.

	\begin{corollary}\label{coro:dfblSubGpToNDflbVal}
		Let $(K,<)$ be an ordered field and let $G := \vnat(K^\times)$. Let furthermore be $H \subseteq G$ an $\Log$-definable convex subgroup such that $\Phi_K^{-1}(H)$ is not $\Lor$-definable. Then $H = \{0\}$ or $H$ is the maximal convex divisible subgroup of $G$.
	\end{corollary}
	\begin{proof}
		By \autoref{theorem:mostDfblSubGpCorrespDfblVal} it follows that $H \subseteq G_0$ since otherwise $\Phi_K^{-1}(H)$ was $\Lor$-definable. Furthermore note that for every prime $p \in \N$ we have $G_0 \subseteq G_p$. Assume $H \not= \{0\}$, else we are done. By \cite[Corollary 4.3]{delon} it follows that there is some prime $p \in \N$ such that $G_p \subseteq H$, because $H$ was $\Log$-definable. But now we have $G_0 \subseteq G_p \subseteq H \subseteq G_0$, hence $H = G_0$. This completes the proof.
	\end{proof}

	Hence there are only two candidates for elements of $\drk_G \setminus \Phi_K(\drk_K)$. To have control of the converse direction, i.e.\ when $\Lor$-definability of a convex valuation implies $\Log$-definability, we now only consider ordered fields with a henselian natural valuation. In this case we can conclusively compare the definable ranks of $K$ and $G$ under the map $\Phi_K$.
	
	\begin{lemma}\label{lemma:henValDrkKEmbedsInDrkG}
		Let $(K,<)$ be an ordered field and let $G = \vnat(K^\times)$. If $\vnat$ is henselian, then $\Phi_K(\drk_K)$ is a final segment of $\drk_G$.
	\end{lemma}
	\begin{proof}
		We first show that for an $\L_{\mathrm{or}}$-definable convex valuation $w$ on $K$ the corresponding convex subgroup $\Phi_K(\OO_w)$ is $\L_{\mathrm{og}}$-definable in $G$. As $w$ is convex, it is a coarsening of $\vnat$. Furthermore, $\OO_w$ is definable in the expansion $(K,<,\OO_\vnat)$. Then, since $\vnat$ is henselian, it follows by \cite[Corollary 4.3 (ii)]{DittmannJahnkeKrappKuhlmann23} that $\Phi_K(\OO_w) = \vnat({\OO_w}^\times)$ is $\L_{\mathrm{og}}$-definable in $G$.
		This establishes that $\Phi_K(\drk_K) \subseteq \drk_G$.
	
		To show that $\Phi_K(\drk_K)$ is a final segment of $\drk_G$, let $H \subseteq G$ be an $\L_{\mathrm{og}}$-definable convex subgroup such that there is an $\L_{\mathrm{or}}$-definable convex valuation $w$ with $\Phi_K(\OO_w) \subseteq H$. We verify that $\Phi_K^{-1}(H)$ is also $\L_{\mathrm{or}}$-definable. By \cite[Proposition 3-1]{DelonLucas1989}, we obtain that $H/\Phi_K(\OO_w)$ is definable in $G/\Phi_K(\OO_w) = w(K^\times)$.
		Therefore \autoref{lem:CoarseningOfDefinableWithDfblSubgp} applies, yielding the desired result.
	\end{proof}

	\begin{theorem}\label{theorem:MainTheorem}
		Let $(K,<)$ be an ordered field such that $\vnat$ is henselian. Denote $G := \vnat(K^\times)$. Then one of the following holds:
		\begin{itemize}
			\item[(1)] $\Phi_K(\drk_K) = \drk_G$ or
		
			\item[(2)] $\Phi_K(\drk_K) = \drk_G \setminus \{\{0\}\}$.
		\end{itemize}
		Furthermore (1) holds if and only if $K\vnat$ is not real closed or there is a prime $p \in \N$ such that $G$ has no non-trivial $p$-divisible convex subgroup.
	\end{theorem}
	\begin{proof}
		By \autoref{lemma:henValDrkKEmbedsInDrkG} it follows that $\Phi_K(\drk_K) \subseteq \drk_G$ is a final segment. Also $K\vnat$ is dense in $\R$ since it is an archimedean real field by the definition of $\vnat$. If $K\vnat$ is not real closed, then it is a proper dense subset of its real closure and in particular not closed in its real closure. By \cite[Theorem 3.1~(3)]{DittmannJahnkeKrappKuhlmann23} this implies that $\vnat$ is definable and therefore $\Phi_K(\drk_K) = \drk_G$ as $\Phi_K(\drk_K)$ is a final segment of $\drk_G$ containing its minimum $\{0\}$. Thus (1) holds in this case.
		
		Assume from now on that $K\vnat$ is real closed. If $p$ is a prime such that $G$ has no non-trivial $p$-divisible convex subgroup, then $G_p = \{0\}$ and by \cite[Proposition 2.6]{delon} it follows that $\vnat$ is definable. This yields (1) as above.
		
		Now also assume that $G$ has a non-trivial $p$-divisible convex subgroup for every prime $p \in \N$. Then $\{0\} \subsetneq G_p$ for all $p$ and $\vnat$ is not definable by \cite[Theorem 4.4]{delon}. Consider $H \in \drk_G \setminus \{\{0\}\}$, if $H \not= G_0$, then by \autoref{coro:dfblSubGpToNDflbVal} we obtain $H \in \Phi_K(\drk_K)$.
		If $G_0 = \{0\}$, then we already obtain (2).
		At last, assume $G_0 \not= \{0\}$ and $H = G_0$. Since $G_p \subseteq G_0$ if and only if $G_p = G_0$ for any prime $p$ and all the $G_p$ are $\Log$-definable, it now follows with \cite[Corollary 4.3]{delon}, that $G_0 \in \drk_G$ if and only if $G_0 = G_p$ for some prime $p$. But then $G_0 \in \Phi_K(\drk_K)$ since $G_0 = G_p = \Phi_K(v_p)$ and $v_p \in \drk_K$ by \cite[Proposition 2.6]{delon}. This now yields (2) and therefore completes the proof.
	\end{proof}

	\begin{corollary}\label{coro:necessaryCondElemEquivToArchimed}
		Let $(k,<)$ be an archimedean ordered field and consider an ordered field $(K,<)$ such that $(K,<) \equiv (k,<)$. If $(K,<)$ is not archimedean, then $K\vnat$ is real closed and $\vnat(K^\times)$ is divisible.
	\end{corollary}
	\begin{proof}
		First note that $\drk_K = \emptyset$ since otherwise there was a formula $\varphi(x;\ul y)$ and parameter $\ul a$ such that $\varphi(x,\ul a)$ defines a non-trivial convex valuation ring on $K$. Then the theory of $K$ would include a sentence stating that there exist parameter $\ul z$ such that $\varphi(x, \ul z)$ defines a non-trivial convex ring. But this sentence is necessarily false on $(k,<)$ as there is no such valuation ring on an archimedean ordered field.
		
		Assume $(K,<)$ is not archimedean, then $\vnat$ is a non-trivial convex valuation on $K$. If $K\vnat$ was not real closed, then it would be definable by \cite[Theorem 3.1.~(iii)]{DittmannJahnkeKrappKuhlmann23}, which would yield a contradiction.
		
		It remains to show that $G := \vnat(K^\times)$ is divisible. From \autoref{coro:dfblSubGpToNDflbVal} we now obtain $\drk_G \subseteq \{ \{0\}, G_0 \}$. We fix a prime $p \in \N$ and assume for contradiction that $G$ is not $p$-divisible. Then $G_0 \subseteq G_p \subsetneq G$ and $G_p$ is an $\Log$-definable proper convex subgroup of $G$ (by \autoref{coro:dfblConvSubgroups}~(4) with $n = p$) and it follows $G_p = G_0$. Furthermore, for any $h \in G \setminus pG$ we obtain that $F_p(h)$ is a definable proper convex subgroup of $G$ with $G_p \subseteq F_p(h)$, so here too holds equality. Following the proof of \autoref{prop:convexSubGKeyLemma}, that means we can choose $H_b = F_p(h) = G_0$ and therefore $G/G_0$ is discrete or not closed in its divisible hull. But then \cite[Theorem 3.1.~(i),(ii)]{DittmannJahnkeKrappKuhlmann23} imply that ${\Phi_K}^{-1}(G_0)$ is $\Lor$-definable and thus ${\Phi_K}^{-1}(G_0) \in \drk_K$. Contradiction.
		
		Hence $G$ must have been $p$-divisible and, since $p$ was an arbitrary prime, it follows that $G$ is divisible.
	\end{proof}

	Note that the converse does not hold, the natural valuation of the non-archi\-medean field $k$ from \autoref{xmpl:drkNonIsomorphic} has a real closed residue field and a divisible value group, but is $\Lor$-definable.

\section{Definability of final segments}\label{section:finalSegments}
	We continue this work with a brief investigation into the definability of final segments of an ordered set $(\Gamma,<)$. One easy to observe fact is the following.
	
	\begin{lemma}\label{lem:princFinSegDfbl}
		Let $\Delta \subseteq \Gamma$ be a final segment such that
		\begin{itemize}
			\item[(a)] $\Delta = \gamma_-$ or
			
			\item[(b)] $\Delta = \gamma_+$
		\end{itemize} for some $\gamma \in \Gamma$.
		Then $\Delta$ is $\L_{\mathrm{<}}$-definable in $\Gamma$.
	\end{lemma}
	
	Next we consider $\Gamma$ to be a densely ordered set or a discretely ordered set. By \cite[Proposition 1.4.~(i),(ii)]{PillaySteinhorn86} those sets are o-minimal, from which we can quickly deduce that all definable final segments are of the form $\gamma_-$ or $\gamma_+$ for some $\gamma \in \Gamma$.
	
	\begin{lemma}\label{lem:denseDiscreteSimple}
		Let $(\Gamma,<)$ be a ordered set and $\Delta \subsetneq \Gamma$ a non-trivial proper definable final segment and assume $\Gamma$ is dense or discrete. Then there is $\gamma \in \Gamma$ such that $\Delta = \gamma_-$ or $\Delta = \gamma_+$.
	\end{lemma}
	\begin{proof}
		By \cite[Proposition 1.4.~(i),(ii)]{PillaySteinhorn86} $\Gamma$ is o-minimal, hence $\Delta$ is a finite union of singletons $s_i$ and open intervals $(l_j,r_j)$. Note that all the intervals must have a lower bound $L-j \in \Gamma$, else $(l_j,r_j)$ is an initial segment of $\Gamma$ and $(l_j,r_j) \subseteq \Delta$ now implies $\Delta = \Gamma$. Consider $s := \min s_i, l := \min l_j$. If $l < s$, then $l = \inf \Delta$ and $\Delta = l_+$, otherwise $s = \min \Delta$ and $\Delta = s_-$.
	\end{proof}
	
	When we consider an ordered set that is neither dense nor discrete, then other final segments can become definable.
	
	\begin{example}\label{xmpl:interestingCut1}
		 Consider the set $\R + \Z$. Here $\emptyset + \Z$ is a final segment without minimum. Furthermore, its complement $\R + \emptyset$ does not have a maximum. But the $\L_{\mathrm{<}}$-formula
		\[ \varphi(x) := \exists y\colon (y<x \wedge \forall z\colon (z \leq y \vee z \geq x)), \]
		i.e.\ $x$ has a predecessor, defines this final segment in $\R + \Z$.
	\end{example}

	In this example, the definable final segment is exactly the maximal convex discrete subordering of $\R + \Z$, i.e.\ the full ordering changes from dense to discrete at this cut. That is no accident as we can conclude from the following result of Rubin.
	
	\begin{fact}\label{fact:Rubin}
		\cite[Corollary 2.3]{Rubin74} Let $\mathcal{A}$ be a model of an expansion $\L$ of $\L_<$ with domain $A$. Let $B \subseteq A$ be convex and $\mathcal{B}$ a substructure of $\mathcal{A}$ with domain $B$. For every $\L$-formula $\phi(v_1,\ldots,v_l,x_1,\ldots,x_k)$ and $\ul a \in (A \setminus B)^k$ there is an $\L$-formula $\phi^*(v_1,\ldots,v_l)$ such that $\mathcal{A} \models \phi(\ul b,\ul a) \Leftrightarrow \mathcal{B} \models \phi^*(\ul b)$ for all $\ul b \in B^l$.
	\end{fact}
	
	\begin{corollary}\label{coro:RubinImplication}
		Let $(\Gamma,<)$ be an ordered set, $\Delta \subseteq \Gamma$ a definable final segment and  $B \subseteq \Gamma$ convex. Then $\Delta \cap B$ is a definable final segment of $(B,<)$.
	\end{corollary}
	\begin{proof}
		Since $\Delta$ is definable there are parameter $\ul a \in \Gamma^n$ and an $\L_<$-formula $\phi(v,\ul x)$ such that $\Delta = \phi(\Gamma, \ul a)$. Without loss of generality we can sort the parameter $\ul a$ such that the first $m$ are from $B$ and the remaining $n-m$ are from $\Gamma \setminus B$. Similarly reorder the variables of $\ul x$ in $\phi(v,\ul x)$ such that they match with the reordered $\ul a$. We can then write $\phi(v,\ul x)$ as $\phi'(v,\ul y,\ul z)$ where $\ul y$ is the first $m$ variables form $\ul x$ and $\ul z$ the remaining $n-m$. Then by \autoref{fact:Rubin} we obtain a formula $\phi^*(v,\ul y)$ such that $b \in \phi'(\Gamma,\ul a)$ if and only if $b \in \phi^*(B,{\ul a}')$ for $b \in B$, where ${\ul a}'$ is the first $m$ parameter from $\ul a$. Hence $\phi^*(B,{\ul a}') = \Delta \cap B$.
	\end{proof}

	With this we can conclude that any definable final segment of $\Gamma$, which is not of the form $\gamma_-$ or $\gamma_+$ for some $\gamma \in \Gamma$, is either disjoint from any convex dense or discrete subset $B$ of $\Gamma$ or fully contain it.
	
	\begin{corollary}\label{coro:localDefblSegm}
		Let $(\Gamma,<)$ be an ordered set, $\Delta \subseteq \Gamma$, $B \subseteq \Gamma$ convex such that $\emptyset \subsetneq \Delta \cap B \subsetneq B$. If $B$ is dense or discrete, then $\Delta$ is of the form $\gamma_-$ or $\gamma_+$ for some $\gamma \in \Gamma$.
	\end{corollary}
	\begin{proof}
		Since $B$ is convex, it follows from \autoref{coro:RubinImplication} that $\Delta \cap B$ is a definable final segment of $B$. By $\emptyset \subsetneq \Delta \cap B \subsetneq B$ it follows that $\Delta \cap B \subsetneq B$ is non-trivial and proper. If $B$ is dense or discrete, then \autoref{lem:denseDiscreteSimple} yields that $\Delta \cap B$  is of the form $\gamma_-$ or $\gamma_+$ for some $\gamma \in B \subseteq \Gamma$. Then so is $\Delta$.
	\end{proof}
	
	That means that through this lens, all interesting definability questions for final segments happen when an ordering is locally neither dense nor discrete. 
	
	\begin{definition}\label{def:Condense}
		(see \cite[Definition 4.1]{Rosenstein82})\ Let $(\Gamma,<)$ be a linear ordering. A \textbf{condensation} of $\Gamma$ is a partition $C$ of $\Gamma$ into convex sets $I$. The condensation becomes a linear ordering with the following relation:
		\[ I_1 < I_2 :\Leftrightarrow \forall a \in I_1, b \in I_2\colon a < b \]
	\end{definition}

	The idea is to define condensations in such a way that exactly points living in a dense interval form singletons in the partition. Then non-dense intervals get `condensed' into single points, through iteration this gradually approximates a dense order. None of the condensations presented in \cite{Rosenstein82} are definable. Instead we propose the following condensation.
	
	\begin{proposition}\label{prop:dfblCondense}
		Let $(\Gamma,<)$ be an ordered set. Consider the $\L_<$-formula $\psi(x;y) := x<y \wedge \forall z\colon \neg(x<z<y)$, i.e.\ $y$ is the successor of $x$. Then the formula
		\[\begin{array}{ccccl}
			\varphi(x,y) := & & x = y & \\
			& \vee & [x<y & \wedge & (\forall \ell\colon x \leq \ell < y \Rightarrow \exists s\colon \psi(\ell,s))\\
			& & & \wedge & (\forall u\colon x < u \leq y \Rightarrow \exists p\colon \psi(p,u))]\\
			& \vee & [y<x & \wedge & (\forall \ell\colon y \leq \ell < x \Rightarrow \exists s\colon \psi(\ell,s))\\
			& & & \wedge & (\forall u\colon y < u \leq x \Rightarrow \exists p\colon \psi(p,u))]
		\end{array}\]
		defines an equivalence relation on $\Gamma$ with convex equivalence classes.
	\end{proposition}
	\begin{proof}
		The formula $\varphi(x,y)$ is certainly reflexive and symmetrical. If $a<b$ and $\Gamma \models \varphi(a,b)$, then for all $c$ with $a<c<b$ we have $\Gamma \models \varphi(a,c)$ and $\Gamma \models \varphi(c,b)$ as the universal quantifiers now only check the conditions for less points.
		
		It only remains to show that if $a< c < b$ and $\Gamma \models \varphi(a,c) \wedge \varphi(c,b)$, then already $\Gamma \models \varphi(a,b)$. From the assumptions we immediately obtain that there are $s,p$ such that $\Gamma \models \psi(a,s) \wedge \psi(p,b)$. If there was a $d$ with $a<d<b$ such that $d$ has no successor or no predecessor, then we would immediately obtain a contradiction to $\varphi(a,c)$ or $\varphi(c,b)$. Hence the relation is transitive.
	\end{proof}
	
	\begin{definition}\label{def:dfblCondense}
		Let $(\Gamma,<)$ be an ordered set. The \textbf{definable condensation} of $\Gamma$ is given by the equivalence classes of the relation defined by $\varphi(x,y)$ from \autoref{prop:dfblCondense}.
	\end{definition}
	
	This condensation is exactly the way we need it. The fibres of points from the condensation are the maximal discrete convex subsets of the original ordering. Therefore all points that are being condensed into a singleton did not have any interesting cuts, since we already treated cuts occurring inside a discrete convex subset with \autoref{coro:localDefblSegm}. Similarly, if this condensation doesn't change a convex subset, then it must already have been dense, leading to the same result.
	Using this condensation, we can now handle more definability questions than before. However, this alone is not enough to describe all definable final segments.
	
	\begin{example}\label{xmpl:condApplication}
		\begin{itemize}
			\item[(i)] Consider the linear ordering $\R + \Z$ from \autoref{xmpl:interestingCut1}. Applying the definable condensation does not change the set on the densely ordered convex subset $\R$, i.e.\ all points there are their own equivalence class. The discrete convex subset $\Z$ on the other hand gets condensed into a single point. We obtain the densely ordered set $\R + \{a\}$. In particular the final segment $\Z$ corresponds to the final segment $\{a\} = a_-$ and is thus definable in the condensation. As a result we obtain here again that this final segment is definable in $\R + \Z$.
			
			\item[(ii)] We now consider the ordered set $\Gamma := \sum_{i \in \R\setminus\{0\}} \Gamma_i$ where $\Gamma_i := \omega$ for $i > 0$ and $\Gamma_j := \omega^*$ for $j < 0$. Then the definable condensation yields $\R \setminus \{0\}$ as result, the final segment $\R^+$ is not definable in this dense ordering. However, the corresponding final segment $\sum_{i \in \R^+} \Gamma_i$ is definable by the formula $\theta(x) := \forall y\colon y > x \Rightarrow [\exists z\colon \psi(y,z)]$ (where $\psi(a,b)$ is the formula from \autoref{prop:dfblCondense} which states that $b$ is the successor of $a$).
		\end{itemize}
	\end{example}
	
	To deal with cases as in \autoref{xmpl:condApplication}~(ii) one could instead consider the definable condensation not just as an ordered set, but as a coloured chain, i.e.\ a linear ordering with additional unary predicates. To do so one can add a label to every point, describing its fibre up to elementary equivalence. Note that the possibilities are limited to, on one hand their cardinality for finite fibres, and on the other hand one of the four models of infinite discrete orderings detailed in \cite[Theorem 2.15.]{robinson}.
	
	However, even then this would not immediately solve all problems, there are linear orderings that do not become dense after any finite iteration of this condensation. What happens in such a case in terms of definability is unclear. A further investigation would greatly exceed the scope of this work and will therefore be left open.

\section{Open questions}
	We conclude with several open questions building on our study. For one, we have seen that in general the definable rank of an ordered field need not agree with that of the value group of its natural valuation, as seen in \autoref{xmpl:ArchG}~(ii). Similarly, the definable rank of an ordered abelian group can be different than that of the value set of its natural valuation. This now leads to the question whether there are relations between the definable ranks on field-, group- and set-level that are always satisfied. More precisely we pose the following question.

	\begin{question}\label{qstn:FieldCondition}
		Given three ordered sets $(A,<), (B,<), (C,<)$. When is it possible to construct an ordered field $K$ such that $\drk_K \cong (A,<)$, $\drk_G \cong (B,<)$ and $\drk_\Gamma \cong (C,<)$?
	\end{question}
	
	We immediately point out that the choice of $(C,<)$ by itself must fulfil very specific conditions. In particular consider the following.
	
	\begin{lemma}
		Let $(\Gamma,<)$ be an ordered set. If there are $\Delta_1, \Delta_2 \in \drk_\Gamma$ with $\Delta_1 \subsetneq \Delta_2$, then there exists $\Delta \in \drk_\Gamma$ with $\Delta_1 \subseteq \Delta \subsetneq \Delta_2$ such that $\Delta$ has a successor.
	\end{lemma}
	\begin{proof}
		Since $\Delta_1 \subsetneq \Delta_2$, there exists some $\gamma \in \Delta_2 \setminus \Delta_1$. Hence $\Delta_1 \subseteq \gamma_+ \subsetneq \gamma_- \subseteq \Delta_2$. Thus, $\gamma_+ \in \drk_\Gamma$ with the successor $\gamma_-$. 
	\end{proof}
	
	That means the definable rank of an ordered set $\Gamma$ has never a dense convex substructure which is not a singleton. We therefore pose a modified version of \autoref{qstn:FieldCondition}.
	
	\begin{question}
		Given an ordered set $(C,<)$. When is it possible to construct an ordered field $K$ such that $\drk_\Gamma \cong (C,<), {\Phi_G}^{-1}(\drk_\Gamma) = \drk_G$ and ${\Phi_K}^{-1}(\drk_G) = \drk_K$ where $G = \vnat(K^\times)$ and $\Gamma = v_G(G \setminus \{0\})$?
	\end{question}
	
	Alternatively one could avoid the question which definable ranks an ordered set may have and instead ask when for a given ordered set $\Gamma$ one can construct an ordered field $K$ with value set $\Gamma$ such that the isomorphisms $\Phi_K, \Phi_G$ preserve definability.
	
	\begin{question}
		Given an ordered set $(\Gamma,<)$. When is it possible to construct an ordered field $K$ such that $v_G(G \setminus \{0\}) = \Gamma$, $\Phi_K(\drk_K) = \drk_G$ and $\Phi_G(\drk_G) = \drk_\Gamma$, where $G = \vnat(K^\times)$?
	\end{question}

	This question has now been answered in \cite[Corollary 4.1]{BoissonneauVogel25}, this is always possible!\\
	
	In \autoref{theorem:mostDfblSubGpCorrespDfblVal} we have established that for any ordered field $(K,<)$ there are at most two definable convex subgroups of its natural value group that do not correspond to a definable convex valuation ring. We have also seen in \autoref{xmpl:drkNonIsomorphic} that there may be definable convex valuation rings that do not correspond to a definable convex subgroup of the natural value group. This naturally gives rise to the final question.
	
	\begin{question}
		Is there a cardinal $\kappa$ such that $|\drk_K \setminus {\Phi_K}^{-1}(\drk_G)| \leq \kappa$ for any ordered field $(K,<)$ and $G = \vnat(K^\times)$?
	\end{question}

\end{document}